\documentclass[oneside, 12pt]{amsart}
\usepackage{amscd, amssymb, amsmath, mathrsfs}
\usepackage[english]{babel}
\usepackage{booktabs}
\usepackage{tikz-cd}
\usepackage[pdftex, colorlinks=true,  citecolor=blue, linkcolor=blue, linktocpage=true]{hyperref}

\newcommand\mylabel[1]{\label{#1}\marginpar{\vspace{-1ex}\medskip\medskip\footnotesize \tt #1}}
\renewcommand\mylabel[1]{\label{#1}}
\newcommand{\mydate}{
\number\day\space
\ifcase\month \or January\or February\or March\or April\or May\or June\or July\or August\or September\or October\or November\or December\fi 
\space\number\year}

\newtheorem{theorem}{Theorem}[section]
\newtheorem*{maintheorem}{Theorem}
\newtheorem{lemma}[theorem]{Lemma}
\newtheorem{proposition}[theorem]{Proposition}
\newtheorem{corollary}[theorem]{Corollary}

\theoremstyle{definition}

\newtheorem*{acknowledgement}{Acknowledgement}

\theoremstyle{remark}

\newcommand{\ZZ}{\mathbb{Z}}

\newcommand{\RR}{\mathbb{R}}
\newcommand{\CC}{\mathbb{C}}
\newcommand{\FF}{\mathbb{F}}

\newcommand{\PP}{\mathbb{P}}
\renewcommand{\AA}{\mathbb{A}}
\newcommand{\GG}{\mathbb{G}}

\newcommand{\ideala}{\mathfrak{a}}
\newcommand{\idealb}{\mathfrak{b}}

\newcommand{\idealg}{\mathfrak{g}}

\newcommand{\shA}{\mathscr{A}}

\newcommand{\shF}{\mathscr{F}}

\newcommand{\Aut}{\operatorname{Aut}}

\newcommand{\can}{\text{\rm  can}}

\newcommand{\id}{{\operatorname{id}}}

\newcommand{\I}{\text{\rm I}}

\newcommand{\Km}{\operatorname{Km}}

\newcommand{\Lie}{\operatorname{Lie}}

\newcommand{\lra}{\longrightarrow}

\newcommand{\MW}{\operatorname{MW}}

\newcommand{\NE}{\operatorname{\overline{NE}}}

\newcommand{\NS}{\operatorname{NS}}

\renewcommand{\O}{\mathscr{O}}

\newcommand{\Pic}{\operatorname{Pic}}

\newcommand{\pr}{\operatorname{pr}}

\newcommand{\quadand}{\quad\text{and}\quad}

\newcommand{\ra}{\rightarrow}

\newcommand{\red}{{\operatorname{red}}}
\newcommand{\Reg}{\operatorname{Reg}}

\newcommand{\Sing}{\operatorname{Sing}}

\newcommand{\Spec}{\operatorname{Spec}}

\newcommand{\pmap}{{[p]}}
\newcommand{\lieg}{\mathfrak{g}}
\newcommand{\liea}{\mathfrak{a}}
\newcommand{\lieb}{\mathfrak{b}}
\newcommand{\liegl}{\mathfrak{gl}}
\newcommand{\ad}{\operatorname{ad}}
\newcommand{\inv}{{-1}}

\begin{document}

\title[ Kummer surfaces with group schemes]
      {Kummer surfaces associated with group schemes}

\author[Shigeyuki Kondo]{Shigeyuki Kond\=o}
\address{Graduate School of Mathematics, Nagoya University, Nagoya 464-8602, Japan}
\curraddr{}
\email{kondo@math.nagoya-u.ac.jp}

\author[Stefan Schr\"oer]{Stefan Schr\"oer}
\address{Mathematisches Institut, Heinrich-Heine-Universit\"at,
40204 D\"usseldorf, Germany}
\curraddr{}
\email{schroeer@math.uni-duesseldorf.de}

% K3 surfaces and Enriques surfaces
% group schemes
% elliptic surfaces

\dedicatory{December 20, 2019}

\begin{abstract}
We introduce Kummer surfaces $X=\Km(C\times C)$ with the group scheme
$G=\mu_2$ acting on the  self-product of the rational cuspidal curve 
in characteristic two. The resulting quotients are normal  surfaces having
a configuration of sixteen rational double points of type $A_1$, together
with a rational double point of type $D_4$. We show that our
Kummer surfaces are precisely the supersingular K3 surfaces with Artin invariant
$\sigma\leq 3$, and characterize them by the existence of a certain configuration of
thirty curves. After contracting suitable curves, they also appear as normal K3-like coverings for simply-connected 
Enriques surfaces.
\end{abstract}

\maketitle
\tableofcontents

%===========================================================
\section*{Introduction}
\mylabel{Introduction}

For each abelian surface $A$ in characteristic $p\neq 2$  the quotient $Z=A/\{\pm 1\}$
by the sign involution is a normal surface with  sixteen rational double points of type $A_1$, and the
minimal resolution of singularities $X=\Km(A)$  is a K3 surface  
called \emph{Kummer surface}. 
Over the complex numbers a K3 surface is a Kummer surface if and only if it contains sixteen disjoint
$(-2)$-curves (Nikulin \cite{Nikulin 1975}).  In this paper we call a non-singular rational curve on a K3 surface
 a \emph{$(-2)$-curve} for simplicity.
If $A=J$ is the jacobian variety of a curve of genus two,
the Kummer surface $X$ contains thirty-two distinguished $(-2)$-curves forming the so-called
$(16_6)$-configuration (e.g.\ Griffiths and Harris \cite{Griffiths; Harris 1978}, Chapter 6, page 787, Figure 21), and the existence of these thirty-two $(-2)$-curves characterizes the Kummer surface associated with a curve of genus 2
(Nikulin \cite{Nikulin 1974}).

Shioda \cite{Shioda 1975} showed that a Kummer surface $X=\Km(A)$ in odd characteristics 
is supersingular if and only if the abelian surface $A$ is supersingular.
For $p=2$, however, 
Shioda \cite{Shioda 1974} and Katsura \cite{Katsura 1978} observed that  
the  singularities on $Z=A/\{\pm 1\}$ are more complicated, and that $X$
is a K3 surface if and only if $A$ is ordinary, and then $\rho(X)\leq 20$.  Indeed, for any supersingular 
abelian surface the quotient acquires an elliptic singularity, and $X$ becomes a rational surface.

The second author \cite{Schroeer 2007} obtained Kummer surfaces in characteristic $p=2$ by replacing the supersingular abelian surface
$A$ by the self-product $C\times C$ of the rational cuspidal curve $C$, which is a \emph{non-normal} surface,
and the constant group $G=\{\pm 1\}$ by the
additive group scheme $G=\alpha_2$, which is \emph{non-reduced}. 
It turns out that for suitable $G$-actions, the quotients $Z=(C\times C)/G$
are normal, and   the minimal resolutions of singularities $X$ are supersingular K3 surfaces
with Artin invariant $\sigma\leq 2$. The configuration of singularities is either $5D_4$  or $D_4+2D_8$.

Recall that for supersingular K3-surfaces,  the isomorphism class of the Picard lattice 
is determined by a single integer $1\leq \sigma\leq 10$ called the \emph{Artin invariant}.
Ogus \cite{Ogus 1979} proved   for odd primes that $\sigma\leq 2$ means that $X$ is Kummer.
For the   Kummer surfaces $X=\Km(C\times C)$ with group scheme $G=\alpha_2$ one also   $\sigma\leq 2$.
Up to isomorphism, there is a unique  supersingular K3 surface with $\sigma=1$.
Dolgachev and the first author \cite{Dolgachev; Kondo 2003} characterized it, among other things, by the existence
of 
%certain genus-one fibrations.
forty-two $(-2)$-curves.
Shimada and Zhang \cite{Shimada; Zhang 2007} showed that every supersingular K3 surface with $\sigma\leq 2$ is isomorphic
to a Kummer surface with $G=\alpha_2$,
by characterizing these Kummer surfaces in terms of the  configurations of twenty-six $(-2)$-curves
with dual graph given in Figure \ref{second curve configuration}.

The main goal of this paper is to extend the construction   $X=\Km(C\times C)$ to the multiplicative group scheme
$G=\mu_2$. It turns out that    $Z=(C\times C)/G$ has only rational double points,
and that their configuration is $16 A_1+D_4$, which is very close to the classical situation over the complex numbers.
The Artin invariant now becomes $\sigma\leq 3$. Indeed, the construction of the Kummer surface  with group scheme $G=\mu_2$
has two moduli coming from the possible embeddings $G\ra \Aut_{C\times C}$. 
Our principal result  is a   characterization  of such Kummer surfaces.

\begin{maintheorem}
(See Thm.\ \ref{configuration and kummer} and \ref{kummer artin})
Let $X$ be a K3 surface in characteristic $p=2$. Then the following are equivalent:
\begin{enumerate}
\item There is an isomorphism $X\simeq \Km(C\times C)$ for 
a Kummer surface  with group scheme $G=\mu_2$.
\item There is  a configuration of thirty $(-2)$-curves on $X$ with simple normal crossings and dual
graph given in \S \ref{Kummer surfaces}, Figure \ref{curve configuration}. 
\item The K3 surface $X$ is supersingular with Artin invariant $\sigma\leq 3$.
\end{enumerate}
\end{maintheorem}

Keum \cite{Keum 1990} showed that  every Kummer surface $X=\Km(A)$ over the complex numbers is the K3-covering
of some Enriques surface $Y$. This was extended to    odd characteristics by  Jang \cite{Jang 2015}.
The first author \cite{Kondo 2018}, \S 3.3 established this also   for our Kummer surface with
Artin invariant $\sigma =1$. Here we extend this to  all Kummer surfaces $X=\Km(C\times C) $
with group scheme $G=\mu_2$: 

\begin{maintheorem} 
(See Thm.\ \ref{k3-like covering}) Let $X$ be a supersingular K3 surface with Artin invariant $\sigma\leq 3$
in  characteristic $p=2$.
Then there is a contraction $X\ra X'$ of twelve $(-2)$-curves such that that the normal K3 surface $X'$
is the K3-like covering of some simply-connected Enriques surface $Y$.
\end{maintheorem}

\medskip
The paper is organized as follows:
In Section \ref{Lie algebras} we recall some facts on group schemes $G$ of height $\leq 1$ and restricted Lie algebras $\lieg$
that will be used throughout. Section \ref{Diagonal actions} contains an analysis of $G$-actions on
the self-product $C\times C$ of the rational cuspidal curve.
This is used in Section \ref{Kummer surfaces} to construct our Kummer surface $X=\Km(C\times C)$ with group scheme $G$,
where we also determine the dual graph for the distinguished curves and compute the Artin invariant.
Our characterization with such configuration of curves occupies Section \ref{Configurations}.
In Section \ref{Artin invariant} we give the  characterization with Artin invariants.

\begin{acknowledgement}
The research of the first author is partially supported by
Grant-in-Aid for Scientific Research (S) No.\ 15H05738.
The research of the second author is conducted in the framework of the research training group
\emph{GRK 2240: Algebro-geometric Methods in Algebra, Arithmetic and Topology}, which is funded
by the DFG. 
\end{acknowledgement}

%===========================================================
\section{Some restricted Lie algebras}
\mylabel{Lie algebras}

We start by discussing some  restricted Lie algebras from a purely algebraic point of view.
Let $k$ be a ground field of characteristic $p>0$. 
Recall that a \emph{restricted Lie algebra} is a Lie algebra $\lieg$,
endowed with an additive self-map $x\mapsto x^\pmap$  called \emph{$p$-map}. The latter  is related to
scalar multiplication,  Lie brackets and vector addition  according to  the following three axioms:
\begin{equation}
\label{restricted lie}
(\lambda x)^\pmap=\lambda^p x,\quad
[x^\pmap,y] = (\ad_x)^p(y),\quad
(x+y)^\pmap = x^\pmap + y^\pmap + \sum_{r=1}^{p-1} s_r(x,y).
\end{equation}
for  all vectors $\ x,y \in \lieg$ and scalars $\lambda\in k$. Here $\ad_x(y)=[x,y]$ is the adjoint representation, and $s_r(x,y)$ are certain universal expressions
involving nested Lie brackets. For example, we have $s_1(x,y)=[x,y]$ in characteristic two,
and $s_1(x,y)= [x,[x,y]]$ in characteristic three. For details, we refer to Demazure and Gabriel \cite{Demazure; Gabriel 1970}, Chapter II, \S7, No.\ 3.

A vector $x\in\lieg$ is called \emph{$p$-closed} if it is nonzero, and $x^\pmap\in \lambda x$ for some scalar $\lambda\in k$.
In other words,  the line  $kx\subset\lieg$ is a restricted Lie subalgebra. 
For each unit $\epsilon\in k^\times$, we get $(\epsilon x)^\pmap=\epsilon^{p-1}\lambda (\epsilon x)$, and we see that
the class of $\lambda$ in $ k/k^{\times p-1}$ depends only on  the line, and not on the vector.
One may regard this class as an ``eigenvalue'' for the $p$-map, and $x\in \lieg$ as an ``eigenvector''.

For each group scheme $G$, the Lie algebra $\lieg=\Lie(G)$ is endowed with a $p$-map in a canonical way.
A group scheme $G$ on which the relative Frobenius map $F:G\ra G^{(p)}$ is trivial
is called \emph{of height $\leq 1$}.
In fact, the functor $G\mapsto \Lie(G)$ is an equivalence between the category of
finite group schemes of height $\leq 1$ and the category of finite-dimensional restricted Lie algebras
(\cite{Demazure; Gabriel 1970}, Chapter II, \S7, No.\ 4).
Such group schemes admit $p$-basis, and their order  $|G|=h^0(\O_G)=p^d$ is given by 
$d=\dim_k(\lieg)$.
In particular, the lines generated by $p$-closed vectors $x\in \lieg$ correspond to subgroup schemes
$H\subset G$ of order $p$. These are twisted form of $H=\mu_p$ or $H=\alpha_p$.
In characteristic $p=2$ we actually have $H=\mu_2$ or $H=\alpha_2$.

We now examine a special type of restricted Lie algebras:
Let $\liea$ be a finite-dimensional restricted Lie algebra, with 
trivial Lie bracket $[a,a']=0$ and trivial $p$-map $a^\pmap = 0$.
Let $\lieb=ke$ be the one-dimensional restricted Lie algebra, with basis vector $e\in\lieb$ and $p$-map
$(\lambda e)^\pmap = \lambda^pe$.
On the vector space sum $\liea\oplus\lieb$,
we define   Lie bracket and $p$-map via the formulas 
\begin{equation}
\label{semidirect formulas}
[a+\lambda e,a'+\lambda'e]=\lambda a'-\lambda'a\quadand (a+\lambda e)^\pmap = \lambda^{p-1}(a+\lambda e).
\end{equation}
The former satisfies the axioms for Lie brackets,
and the resulting Lie algebra is written as $\lieg=\liea\rtimes\lieb$. This is indeed the semi-direct product
with respect to the homomorphism of Lie algebras $\lieb\ra\liegl(\liea)$, $e\mapsto\id_\liea$.

\begin{proposition}
\mylabel{semidirect product}
The above $p$-map endows  $\lieg=\liea\rtimes\lieb$ with the structure of a restricted Lie algebra,
such that the inclusions of $\liea$ and $\lieb$ are homomorphisms of restricted Lie algebras.
Moreover, this $p$-map is unique, and each vector in $\lieg$ is $p$-closed.
\end{proposition}

\proof
Uniqueness follows from the axioms \eqref{restricted lie}. Since the homomorphism of Lie algebras
$\rho:\idealb\ra\liegl(\liea)$
satisfies $ \rho(x^\pmap)=\rho(x)^p$, the existence of a $p$-map 
follows from Strade and Farnsteiner \cite{Strade; Farnsteiner 1988}, Theorem 2.5.
To check that it is given by our formula, it suffices to treat the case $\ideala=ka$, which was
verified in loc.\ cit.\ example 4 on page 72.
From \eqref{semidirect formulas} we see that each vector in the semidirect product is $p$-closed.
\qed

\medskip
Let $G$ be the finite group scheme of height $\leq 1$ with $\Lie(G)=\idealg$. Then 
the closed subschemes $H\subset G$ of order $p$ correspond to lines in the vector space
$\lieg$, or points on the projectivization $\PP(\lieg)$.

Now consider the direct sum $\lieg\oplus\lieg$ of restricted Lie algebras.
Here the Lie bracket and $p$-map are given by
$$
[(x,x'),(y,y')]=([x,y],[x',y'])\quadand (x,x')^\pmap = (x^\pmap,x'^\pmap).
$$
A straight-forward computation shows:

\begin{proposition}
\mylabel{p-closed on sum}
A vector $(x,x')\in\lieg\oplus\lieg$ whose entries are non-zero is $p$-closed if and only if $x=a+\lambda e$
and $x'=a'+\lambda e$ for some common scalar $\lambda\in k$. In this case,
we have $(x,x')^\pmap = \lambda (x,x')$.
\end{proposition}

We see that the set of  $p$-closed vectors is the union of the three restricted Lie subalgebras
$\lieg\oplus 0$ and  $0\oplus \lieg$  and $(\liea\oplus\liea)\rtimes\lieb$,
where the latter is formed from the diagonal embedding $\lieb\subset\lieb\oplus\lieb$.

%===========================================================
\section{Diagonal actions and rational points}
\mylabel{Diagonal actions}

Let $k$ be a ground field of characteristic $p=2$, and consider
the rational cuspidal curve
$$
C=\Spec k[u^2,u^3]\cup \Spec k[u^{-1}].
$$
As explained by the second author in  \cite{Schroeer 2007}, Section 3 the sheaf $\Omega^1_{C/k}$ of K\"ahler differentials
is invertible modulo torsion, and the dual sheaf $\Theta_{C/k}$ is
invertible of degree four. By Riemann--Roch,    $H^0(C,\Theta_{C/k})$
is four-dimensional. As a  restricted Lie algebra this is  a semidirect product $\liea\rtimes\lieb$
studied in the previous section,
where the first factor $\liea$ is generated by the vector fields $u^{-2}D_u,  D_u, u^2D_u$
and the second factor $\lieb$ is generated by $uD_u$.
Here the derivation $D_u:\O_C\ra\O_C$ is determined by $D_u(u)=1$. Note that the basis vectors can be rewritten as
$$
u^{-2}D_u=u^{-4}D_{u^\inv},\quad
D_u=u^{-2}D_{u^\inv},\quad 
u^2D_u=D_{u^\inv} \quadand 
uD_u=u^{-1}D_{u^\inv}.
$$
By  Proposition \ref{semidirect product}, each non-zero vector field $\delta\in H^0(C,\Theta_{C/k})$ is $p$-closed, and  thus
defines a faithful action of the height-one group scheme $G$  with $\Lie(G)=k\delta$.
Note that we have $G\simeq\alpha_2$ if $\delta\in\ideala$ and $G\simeq\mu_2$ else.

We now consider the self-product $C\times C$, which is a non-normal integral surface. 
As discussed in Section \ref{Lie algebras}, the restricted Lie algebra  
\begin{equation}
\label{decomposition}
H^0(C\times C,\Theta_{C\times C/k}) =(\liea\rtimes\lieb)\oplus (\liea\rtimes\lieb)
\end{equation}
contains the restricted Lie subalgebra $(\liea\oplus\liea)\rtimes\lieb$, whose elements
have the form
\begin{equation}
\label{delta}
\delta=(\lambda_4 u^{-4} +\lambda_2u^{-2} + \lambda_0)D_{u^\inv} + (\mu_4 v^{-4} +\mu_2v^{-2} + \mu_0)D_{v^\inv} + \tau(uD_u+vD_v)
\end{equation}
for some scalars $\lambda_i,\mu_i$ and a common $\tau=\lambda_1=\mu_1$.
Here we use the two indeterminates $u,v$ to describe the first and second factor in $C\times C$, respectively.

Now let $G$ be a local group scheme order $p$ acting faithfully on $C\times C$.
In light of the decomposition \eqref{decomposition}, this is the diagonal action stemming from $G$-actions on the factors.
According to Proposition \ref{p-closed on sum},
such actions arise from the non-zero vector field $\delta$ as in \eqref{delta}.

Consider the quotient $Z=(C\times C)/G$, which is an integral   surface. 
The second projection $\pr_2:C\times C\ra C$ induces a morphism
$$
Z=(C\times C)/G\lra C/G=\PP^1.
$$
Here the projective line is given by $\PP^1=\Spec k[v^2]\cup \Spec k[v^{-2}]$. Write $K=k(v^{-2})$
for its function field.
The generic fiber $Z_K= Z\otimes_{\O_{\PP^1}}K$ is a twisted form of the rational cuspidal curve $C_K=C\otimes_kK$.
We compute its $K$-rational points in dependence of the vector field $\delta$:

\begin{proposition}
\mylabel{system of equations}
The set of $K$-rational points in the regular locus $\Reg(Z_K)$ 
corresponds to the solution $(\alpha,\beta)\in K^2$ of the system of   equations
$$
\lambda_4\alpha^4+\mu_4\alpha=0,\quad
\lambda_2\alpha^2+\mu_2\alpha=0, \quad
\lambda_4\beta^4+\lambda_2\beta^2+\lambda_0  + \tau\beta=\mu_0 \alpha.
$$
\end{proposition}

\proof 
Write  $L=k(v^{-1})$ for the function field of the rational cuspidal curve $C=\Spec k[v^2,v^3]\cup\Spec k[v^{-1}]$.
Each $K$-rational point on the regular locus of $Z_K\subset Z$ has as preimage on   $C_L\subset C\times C$ a $G$-stable $L$-valued point.
As explained by the second author in  \cite{Schroeer 2007}, proof for Proposition 7.2, 
these correspond to  $G$-equivariant  $K$-morphisms 
$\Spec(L)\ra \Reg(C_K)=\Spec K[u^\inv]$. To make the actions explicit we consider the polynomials
\begin{gather*}
P(u^\inv)= \lambda_4 u^{-4} +\lambda_2u^{-2} + \lambda_0+\tau u^\inv,\\
Q(v^\inv)= \mu_4 v^{-4} +\mu_2v^{-2} + \mu_0 + \tau v^\inv
\end{gather*}
that appear in the vector field \eqref{delta}.
Then $G$ acts on $C_K$ via the vector field $P(u^\inv)D_{u^\inv}$ and on $\Spec(L)$ by  the derivation $Q(v^\inv)D_{v^\inv}$.
The morphism $\Spec(L)\ra \Reg(C_K)$ is given by a homomorphism of $K$-algebras
$$
K[u^{-1}]\lra L,\quad u^{-1}\longmapsto \alpha v^{-1}+\beta.
$$
for some   $\alpha,\beta\in K$.
The $G$-equivariance means that the substitution $u^\inv=\alpha v^{-1}+\beta$ turns the derivation $P(u^\inv)D_{u^\inv}$
into the derivation $Q(v^\inv)D_{v^\inv}$.
The latter condition boils down to the equation  $P(\alpha v^\inv+\beta)=\alpha Q(v^\inv)$.
Comparing coefficients in this polynomial equation for the indeterminate $v^\inv$ gives the desired system of   equations
in $\alpha,\beta$.
\qed

\medskip
If $\lambda_4\neq 0$, one easily sees that the number of solutions $\alpha\in k$
for the first two equations $ \lambda_4\alpha^4+\mu_4\alpha=0$ and $\lambda_2\alpha^2+\mu_2\alpha=0$
is of the form $2^m$ for some integer $0\leq m\leq 2$.
This leads to a formula for the number of rational points on the generic fiber:

\begin{corollary}
\mylabel{rational points}
Suppose that $\lambda_4,\mu_4\neq 0$, that the generic fiber
 $Z_K$ is normal,  and that $k$ is algebraically closed. Let $0\leq m\leq 2$ be as above.
Then the number of rational points
in the generic fiber for the morphism $Z\ra\PP^1$ is given by 
$$
|Z_K(K)| = \begin{cases}
2^{2+m}	& \text{if $G=\mu_2$;}\\
2^{1+m}	& \text{if $G=\alpha_2$ and $\lambda_2\neq 0$;}\\
2^m	& \text{if $G=\alpha_2$ and $\lambda_2=0$.}
\end{cases}
$$
Moreover, in every case each $0\leq m\leq 2$   occurs for a suitable choice of $\lambda_2,\mu_2\in k$.
\end{corollary}

\proof
First note that   all the solutions $\alpha,\beta\in K$
for the equations in Proposition \ref{system of equations}
already lie in $k$, because this field is relatively algebraically
closed in $K$. Moreover, the set   $Z_K(K)$ is contained in  $\Reg(Z_K)$,
because the curve $Z_K$ is normal.

Suppose first that $G=\mu_2$, in other words $\tau\neq 0$. Then the third  equation 
$\lambda_4\beta^4+\lambda_2\beta^2+\lambda_0  + \tau\beta=\mu_0 \alpha$ is separable, thus   for each solution    $\alpha$
of the first two equation one gets 
four solutions   $\beta$ of the third equation. This gives the desired formula $|Z_K(K)|=  2^m\cdot 4$.
For $\lambda_2$ generic the non-zero solution $\alpha=\mu_2/\lambda_2$ of the second equation is not a solution of the first equation,
thus  $m=0$. The other extreme $\lambda_2=\mu_2=0$ yields
$m=2$, and for suitable choices of $\lambda_2,\mu_2$  we also get $m=1$.
The assertion for $G=\alpha_2$  where $\tau=0$ is similar  and left to the reader. 
\qed

%===========================================================
\section{Kummer surfaces associated with group schemes}
\mylabel{Kummer surfaces}

We keep the assumptions of the previous section, such that the group scheme $G$ acts on the self-product $C\times C$
of the rational cuspidal curve in characteristic $p=2$. The action is  given by some global  vector field
\begin{equation}
\label{derivation}
\delta=(\lambda_4 u^{-4} +\lambda_2u^{-2} + \lambda_0)D_{u^\inv} + (\mu_4 v^{-4} +\mu_2v^{-2} + \mu_0)D_{v^\inv} + \tau(uD_u+vD_v).
\end{equation}
with seven coefficients $\lambda_i,\mu_i,\tau\in k$.  We assume that
the action faithful on each factor.
For the sake of exposition, we also assume that $k$ is algebraically closed.
Write $Z=(C\times C)/G$ for the resulting quotient, which is an integral surface.
As in   \cite{Schroeer 2007}, Proposition 4.3 one verifies:

\begin{proposition}
The following conditions are equivalent:
\begin{enumerate}
\item The  integral scheme $Z=(C\times C)/G$ is normal
\item Both coefficients $\lambda_4,\mu_4\in k$ are non-zero. 
\item The   $G$-actions on the two factors $C$ are free at the singular point.
\end{enumerate}
\end{proposition}
 
From now on we assume that indeed $\lambda_4,\mu_4\neq 0$, and  proceed to study the resulting normal surface 
$Z=(C\times C)/G$. 
The \emph{fixed scheme} $(C\times C)^G$ for the group scheme action thus lies in the regular locus, 
and is thus the zero-scheme for the equations $\lambda_4 u^{-4} +\lambda_2u^{-2} + \lambda_0 +\tau u^\inv=0$
and $\mu_4 u^{-4} +\mu_2u^{-2} + \mu_0 +\tau u^\inv=0$. This is a finite subscheme of length $l=16$.
The point   on $C\times C$ that lies over the singular point  in both factors   is given by $u^2=u^3=v^2=v^3=0$
and called the \emph{quadrupel point}. 

\begin{proposition}
The $G$-action on $C\times C$ and the normal surface $Z=(C\times C)/G$ has the following properties:
\begin{enumerate}
\item
The singular locus $\Sing(Z)$ consists
of the images of the fixed points, together with the image of the quadruple point. 
The latter is always a rational double point of type $D_4$.
\item
For $G=\mu_2$, the fixed scheme $(C\times C)^G$ is reduced, and its image on $Z$ consists of $l=16$ points,
which are rational double points of type $A_1$.
\item 
For $G=\alpha_2$ the fixed scheme is non-reduced, and its image on $Z$ either consists of four rational double points of type $D_4$,
or two rational double points of type $D_8$, or one elliptic singularity. The latter holds if and only if $\lambda_2=\mu_2=0$.
\item
The minimal resolution of singularities $X\ra (C\times C)/G$ is a K3 surface if and only if 
all singularities on $(C\times C)/G$ are rational, that is,  $\lambda_2,\mu_2\in k$ do not vanish simultaneously.  
\end{enumerate}
\end{proposition}

\proof
The case  that $G$ is additive is already treated by the second author in \cite{Schroeer 2007}, Section 5 and 6, under the assumption
$\lambda_4=\mu_4=1$ and $\lambda_0=\mu_0$. The general case works in virtually the same; it can also be
reduced to the special case by using a coordinate change.

\newcommand{\tP}{\tilde{P}}
\newcommand{\tQ}{\tilde{Q}}
The case  that $G$ is multiplicative is  completely analogous. Let us   make the
singularities arising from the fixed points explicit:  The fixed points
are given by the vanishing of the separable polynomials $P=\lambda_4 u^{-4} +\lambda_2u^{-2} + \lambda_0+\tau u^{-1}$
and $Q=\mu_4 v^{-4} +\mu_2v^{-2} + \mu_0+\tau v^{-1}$. Moreover, 
the kernel  $R\subset k[u^{-1},v^{-1}]$  of the derivation $\delta$ is generated as $k$-algebra
by the elements $a=u^2$ and $b=v^2$ and 
$$
c=\tau u^\inv v^\inv + (\mu_4 v^{-4} +\mu_2v^{-2} + \mu_0) u^\inv + (\lambda_4 u^{-4} +\lambda_2u^{-2} + \lambda_0) v^\inv.
$$
These three generators are subject to the single relation
$$
c^2 +  \tau^2ab + (\mu_4^2 b^4 +\mu_2^2b^2 + \mu_0^2) a + (\lambda_4^2 a^4 +\lambda_2^2a^2 + \lambda_0^2)b=0.
$$
One easily checks that this equation indeed defines   sixteen rational double points of type $A_1$.
\qed

\medskip 
From now on, we assume that the coefficients $\lambda_2,\mu_2\in k$ do not vanish simultaneous,
such that the minimal resolution $X\ra (C\times C)/G$ is a supersingular K3 surface. 
Note that the group scheme $G$, the  action on $C\times C$ and   the resulting
 K3 surface $X$
all depend  on the   vector field \eqref{derivation}.  By abuse of notation  
we write $$X=\Km(C\times C)$$ and call it the 
\emph{Kummer surface associated with group scheme $G$}. Note that we either have $G=\mu_2$ or $G=\alpha_2$.

The vector field in \eqref{derivation} depends on the seven   parameters $\lambda_i,\mu_i,\tau\in k$. However,
the isomorphism class of the
Kummer surface $X=\Km(C\times C)$ has only   two   moduli, because the image of the embedding  $G\ra \Aut_{C\times C}$
depends only on the line $k\delta\subset H^0(C\times C,\Theta_{C\times C})$ rather than the vector $\delta$,
and the canonical action of $\GG_a\rtimes\GG_m$ on the affine line $\AA^1$, which extends to an action on $C$,
yields   re-parameterization
of the indeterminates $u^\inv$ and $v^\inv$.

Now suppose  $X=\Km(C\times C)$ is  a Kummer surface with group scheme $G=\mu_2$. Owing to its construction
as minimal resolution of singularities $X\ra (C\times C)/G$, this K3 surface contains thirty \emph{distinguished curves}
\begin{equation}
\label{thirty curves}
E_{ij}, \quad E_r,\quad C_s,\quad C_s' \qquad (1\leq i,j\leq 4,\quad 0\leq r\leq 3, \quad 0\leq s\leq 4)
\end{equation}
defined as follows:
The two projections $\pr_1,\pr_2:C\times C\ra C$ induce two fibrations $f,f':X\ra \PP^1$.
The   $E_{ij}\subset X$ are the exceptional curves lying over the images $(a_i,a_j')\in\PP^1\times\PP^1$
of the sixteen fixed points in $C\times C$. The $E_r\subset X$ are  the  exceptional curves over the image of the quadrupel point.
Finally, the    $C_s,C_s'\subset X$ are the strict transforms of the fibers $f^{-1}(a_s)$ and $f'^{-1}(a_s)$, respectively.

For any genus-one fibration $f:S\to \PP^1$ on a K3 surface $S$ with a section 
$O\subset S$,
the \emph{trivial lattice} $T(S/\PP^1)$ is the sublattice inside $\Pic(S)$ generated by the irreducible 
components of the closed fibers, together with the chosen section.

\begin{figure}
{\small
\begin{tikzpicture}
[node distance=2cm, label distance=-.1cm]
%\draw[help lines] (0,0) grid (12,12);
\tikzstyle{vertex}=[circle, draw, inner sep=0mm, minimum size=2ex]
\tikzstyle{vertexb}=[circle, draw, inner sep=0mm, minimum size=2ex, fill=black]

\node[vertexb]	(C1)  	at (1,6) 		[label=left:{$C_1$}]{};
\node[vertex]	(E11)  	at (2,10) 		[label=  below right:{$E_{11}$}]{};
\node[vertex]	(E12)	[below of=E11]	[label={ below right, label distance=-.1cm :{$E_{12}$}}]{};
\node[vertex]	(E13)	at (2,4)		[label=above right:{$E_{13}$}]{};
\node[vertex]	(E14)	[below of=E13]	[label=above right:{$E_{14}$}]{};

\node[vertexb]	(C2)  	at (3,6) 		[label=left:{$C_2$}]{};
\node[vertex]	(E21)  	at (4,10) 		[label=below right:{$E_{21}$}]{};
\node[vertex]	(E22)	[below of=E21]	[label=below right:{$E_{22}$}]{};
\node[vertex]	(E23)	at (4,4)		[label=above right:{$E_{23}$}]{};
\node[vertex]	(E24)	[below of=E23]	[label=above right:{$E_{24}$}]{};

\node[vertexb]	(C3)  	at (9,6) 		[label=right:{$C_3$ }]{};
\node[vertex]	(E31)  	at (8,10) 		[label=below left:{$E_{31}$}]{};
\node[vertex]	(E32)	[below of=E31]	[label=below left:{$E_{32}$}]{};
\node[vertex]	(E33)	at (8,4)		[label=above left:{$E_{33}$}]{};
\node[vertex]	(E34)	[below of=E33]	[label=above left:{$E_{34}$}]{};

\node[vertexb]	(C4)  	at (11,6) 		[label=right:{$C_4$}]{};
\node[vertex]	(E41)  	at (10,10) 		[label=below left:{ $E_{41}$}]{};
\node[vertex]	(E42)	[below of=E41]	[label=below left:{$E_{42}$}]{};
\node[vertex]	(E43)	at (10,4)		[label=above left:{$E_{43}$}]{};
\node[vertex]	(E44)	[below of=E43]	[label=above left:{$E_{44}$}]{};

\node[vertex]	(E0)  	at (6.5,5.5)	[label=above left:{$E_0$}]{};
\node[vertex]	(E1)  	at (6,5) 		[label=left:{$E_1$}]{};
\node[vertex]	(E2)  	at (7,6) 		[label=above:{$E_2$}]{};
\node[vertex]	(E3)  	at (7,5)		[label=above right:{$E_3$}]{};

\node[vertexb]	(C'1)  	at (6,11) 		[label=above:{$C'_1$}]{};
\node[vertexb]	(C'2)  	at (6,9) 		[label=above:{$C'_2$}]{};
\node[vertexb]	(C'3)  	at (6,3) 		[label=below:{$C'_3$}]{};
\node[vertexb]	(C'4)  	at (6,1) 		[label=below:{$C'_4$}]{};

\node[vertexb]	(C0)  	at (6,7) 		[label=above:{$C'_0$}]{};
\node[vertexb]	(C'0)  	at (5,6) 		[label=left:{$C_0$}]{};

\draw (E11)--(C1)--(E12)--(C1)--(E13)--(C1)--(E14);
\draw (E21)--(C2)--(E22)--(C2)--(E23)--(C2)--(E24);
\draw (E31)--(C3)--(E32)--(C3)--(E33)--(C3)--(E34);
\draw (E41)--(C4)--(E42)--(C4)--(E43)--(C4)--(E44);

\draw (E11)--(C'1)--(E21)--(C'1)--(E31)--(C'1)--(E41);
\draw (E12)--(C'2)--(E22)--(C'2)--(E32)--(C'2)--(E42);
\draw (E13)--(C'3)--(E23)--(C'3)--(E33)--(C'3)--(E43);
\draw (E14)--(C'4)--(E24)--(C'4)--(E34)--(C'4)--(E44);

\draw (C1)--(C0)--(C2)--(C0)--(C3)--(C0)--(C4)--(C0)--(E0);
\draw (C'1)--(C'0)--(C'2)--(C'0)--(C'3)--(C'0)--(C'4)--(C'0)--(E0);
\draw (E1)--(E0)--(E2)--(E0)--(E3);

\end{tikzpicture}
}
\caption{Dual graph for the thirty  distinguished $(-2)$-curves on Kummer surface  $\Km(C\times C)$ associated with group scheme $G=\mu_2$}
\label{curve configuration}
\end{figure}
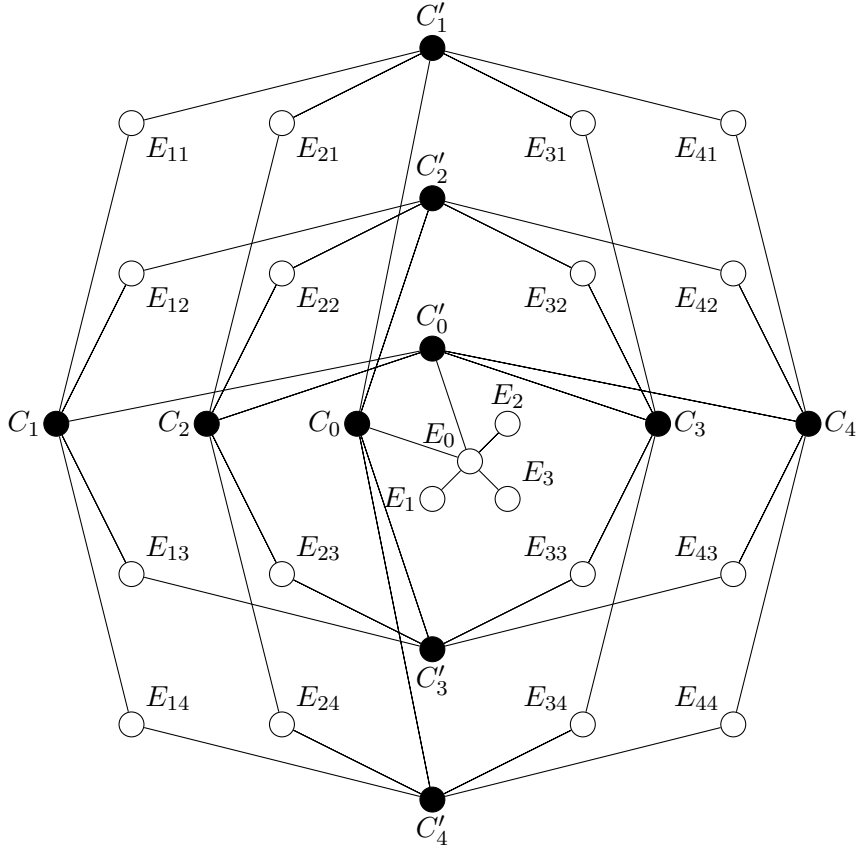

\begin{proposition}
\mylabel{distinguished curves}
For the group scheme $G=\mu_2$, the 
thirty distinguished curves on the Kummer surface $X=\Km(C\times C)$ listed in \eqref{thirty curves} 
form a configuration of $(-2)$-curves with simple normal crossings
whose dual graph is depicted in Figure \ref{curve configuration}.
Moreover, the  Kummer  surface $X$ has Picard number $\rho=22$.
\end{proposition}

\proof
This can be checked with  a local computation with    rings of invariants. 
However, we can also argue that subconfigurations
like $C_1+E_{11}+\ldots+E_{14}$ appear as  set-theoretical fibers for genus-one fibrations $f:X\ra \PP^1$.
Hence the components must be $(-2)$-curves with simple normal crossings, and the shape of the dual graph follows.

The above genus-one fibration $f:X\ra\PP^1$ is actually induced by the first projection $\pr_1:C\times C\ra C$, hence it is quasielliptic.
It has five reducible fibers with Kodaira symbol $\I_0^*$, and $C'_1\subset X$ is a section.
In turn, the trivial lattice $T(X/\PP^1)\subset\Pic(X)$, which is generated by the vertical curves disjoint from $C'_1$, together
with a fiber and $C'_1$,
has rank $22=5\cdot 4+2$. It follows that our K3 surface $X$ has Picard number $\rho=22$.
\qed

% \medskip
% K3-surfaces with Picard number $\rho=22$ are called \emph{supersingular}.
% We just observed that our Kummer surface $X=\Km(C\times C)$ has  a jacobian 
% quasielliptic fibration $f:X\ra\PP^1$ having five reducible fibers with
% Kodaira symbol $\I_0^*$, and the Mordell--Weil group is an elementary abelian $2$-group of rank $r\geq 2$.
% According to \cite{Dolgachev; Kondo 2003}, Theorem 1.1 this must be the unique supersingular K3 surface with $\sigma=1$,
% provided that $r=4$.

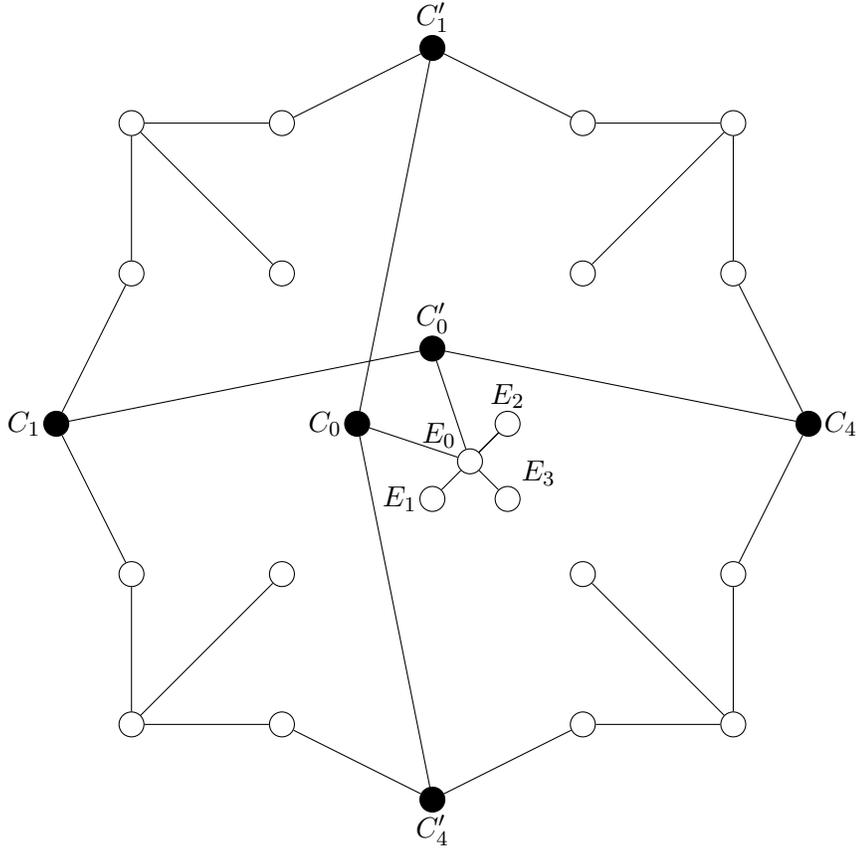
\begin{figure}
{\small
\begin{tikzpicture}
[node distance=2cm, label distance=-.1cm]
%\draw[help lines] (0,0) grid (12,12);
\tikzstyle{vertex}=[circle, draw, inner sep=0mm, minimum size=2ex]
\tikzstyle{vertexb}=[circle, draw, inner sep=0mm, minimum size=2ex, fill=black]

\node[vertexb]	(C1)  	at (1,6) 		[label=left:{$C_1$}]{};
\node[vertex]	(E11)  	at (2,10) 		[label=  below right:{}]{};
\node[vertex]	(E12)	[below of=E11]	[label={ below right, label distance=-.1cm :{}}]{};
\node[vertex]	(E13)	at (2,4)		[label=above right:{}]{};
\node[vertex]	(E14)	[below of=E13]	[label=above right:{}]{};

\node[vertex]	(E21)  	at (4,10) 		[label=below right:{}]{};
\node[vertex]	(E22)	[below of=E21]	[label=below right:{}]{};
\node[vertex]	(E23)	at (4,4)		[label=above right:{}]{};
\node[vertex]	(E24)	[below of=E23]	[label=above right:{}]{};

\node[vertex]	(E31)  	at (8,10) 		[label=below left:{}]{};
\node[vertex]	(E32)	[below of=E31]	[label=below left:{}]{};
\node[vertex]	(E33)	at (8,4)		[label=above left:{}]{};
\node[vertex]	(E34)	[below of=E33]	[label=above left:{}]{};

\node[vertexb]	(C4)  	at (11,6) 		[label=right:{$C_4$}]{};
\node[vertex]	(E41)  	at (10,10) 		[label=below left:{}]{};
\node[vertex]	(E42)	[below of=E41]	[label=below left:{}]{};
\node[vertex]	(E43)	at (10,4)		[label=above left:{}]{};
\node[vertex]	(E44)	[below of=E43]	[label=above left:{}]{};

\node[vertex]	(E0)  	at (6.5,5.5)	[label=above left:{$E_0$}]{};
\node[vertex]	(E1)  	at (6,5) 		[label=left:{$E_1$}]{};
\node[vertex]	(E2)  	at (7,6) 		[label=above:{$E_2$}]{};
\node[vertex]	(E3)  	at (7,5)		[label=above right:{$E_3$}]{};

\node[vertexb]	(C'1)  	at (6,11) 		[label=above:{$C'_1$}]{};
\node[vertexb]	(C'4)  	at (6,1) 		[label=below:{$C'_4$}]{};

\node[vertexb]	(C'0)  	at (6,7) 		[label=above:{$C'_0$}]{};
\node[vertexb]	(C0)  	at (5,6) 		[label=left:{$C_0$}]{};

\draw (E11)--(E12)--(C1)--(E13)--(E14);
\draw (C'1)--(C0)--(C'4);
\draw (C1)--(C'0)--(C4);
\draw (C0)--(E0)--(C'0);
\draw (E41)--(E42)--(C4)--(E43)--(E44);

\draw (E11)--(E21)--(C'1)--(E31)--(E41);
\draw (E11)--(E22);
\draw (E41)--(E32);
\draw (E44)--(E33);
\draw (E14)--(E23);
\draw (E14)--(E24)--(C'4)--(E34)--(E44);

\draw (E1)--(E0)--(E2)--(E0)--(E3);
\end{tikzpicture}
}
\caption{Dual graph for the twenty-six distinguished $(-2)$-curves on  $\Km(C\times C)$ associated with group scheme $G=\alpha_2$ and generic action}
\label{second curve configuration}
\end{figure}

\medskip
Note that removing the curves $C_0, C_0', E_0,\ldots,E_3$ yields the dual graph that occurs in the classical Kummer surface
$X=\Km(E_1\times E_2)$ in characteristic $0$ attached to the product of ellliptic curves, where the quotients
$Z=(E_1\times E_2)/\{\pm1\}$ acquires sixteen rational double points of type $A_1$.
For the sake of completeness,
we depict the distinguished curves on the Kummer surface 
$X=\Km(C\times C)$ associated with group scheme $G=\alpha_2$ and generic action  in 
Figure \ref{second curve configuration}.

\newcommand{\disc}{\operatorname{disc}}
\newcommand{\sign}{\operatorname{sign}}
A finitely generated free abelian group   endowed with a non-degenerate $\ZZ$-valued bilinear form  is called
a \emph{lattice}. 
A lattice $L$ is \emph{$p$-elementary} if the discriminant group $L^*/L$ is an elementary abelian $p$-group.  Here $L^* = {\rm Hom}(L,\ZZ)$.
Let $L$ be a $2$-elementary even lattice of rank $r\geq 2$ with signature $\sign(L)=(1,r-1)$, and 
assume that the discriminant bilinear form $b_L:L^*/L\times L^*/L\ra \frac{1}{2}\ZZ/\ZZ$ 
is alternating,
in other words, the discriminant quadratic form $q_L:L^*/L\ra\frac{1}{2}\ZZ/2\ZZ$ factors over $\ZZ/2\ZZ$.
Regarding the latter as a quadratic form with values in $\ZZ/2\ZZ=\frac{1}{2}\ZZ/\ZZ$,
we see that $b_L$ becomes its associated bilinear form.
Thus $q_L$ is non-degenerate and   the $\FF_2$-vector space $L^*/L$ is symplectic and thus even-dimensional, so  
$\disc(L)=(-1)^{r-1}2^{2\sigma}$ for some integer $\sigma\geq 0$.

According to Nikulin (\cite{Nikulin 1980}, Theorem 3.6.2),  the even indefinite $2$-elementary lattice $L$ is determined up to isometry
by the numerical invariants $r\geq 2$ and $\sigma\geq 0$.
Moreover,  such a lattice with invariants $r\geq 2$ and $\sigma\geq 0$ exists if and only 
if $2\sigma\leq r$ and  furthermore $r\equiv 2$ modulo $4$, where this congruence already holds modulo $8$ in the boundary
cases $\sigma=0$ and $2\sigma=r$.
The discriminant quadratic form $q_L:L^*/L\ra \ZZ/2\ZZ$ is   given by  
\begin{equation}
\label{discriminat quadratic form}
q_L=\begin{cases}
% q_{\sigma}		& \text{if $r\equiv 2$ mod  $8$};\\
% q_{\sigma-1}+t_1^2+t_1t_2+t_2^2	& \text{if $r\equiv 6$ mod  $8$}
q(t_1,\ldots,t_{2\sigma}) 				& \text{if $r\equiv 2$ mod  $8$};\\
q(t_1,\ldots,t_{2\sigma})+t_1^2+t_2^2	& \text{if $r\equiv 6$ mod  $8$}
\end{cases}
\end{equation}
with the polynomial $q(t_1,\ldots,t_{2\sigma}) = t_1t_2+\ldots+t_{2\sigma-1}t_{2\sigma}$.   
Recall that up to isometry, there are two non-degenerate quadratic forms on $V=\FF_2^{\oplus 2\sigma}$ with $\sigma\geq 1$,
which are   distinguished by the \emph{Arf invariant}, or equivalently by the \emph{Witt index}.

An even overlattice $L\subset L'$ corresponds to a subgroup $H\subset L^*/L$ on which 
the discriminant quadratic form $q_L$ vanishes, 
via the assignment $H=L'/L$. Such $H$  is called \emph{totally  singular}  
(Bourbaki's terminology \cite{A 9}, \S9, No.\ 2, Definition 2).
The discriminant group for the overlattice $L'$ is given by  the subquotient $H^\perp/H$, so the bilinear form $b_{L'}$ stays
alternating. In turn, the invariants for the overlattice $L'$ are given by the formulas
$r'=r$ and $\sigma'=\sigma-\dim_{\FF_2}(H)$. It is not difficult to count the number of overlattices:

\begin{lemma}
\mylabel{overlattices}
The number $n\geq 0$ of overlattices $L\subset L'$ with index $[L':L]=2$ is given by the formula
$n=2^{2\sigma-1} + (-1)^{\epsilon} 2^{\sigma-1}-1$ for the exponent $\epsilon=(r-2)/4$. 
\end{lemma}

\proof
The number of totally isotropic subgroups $H\subset L^*/L$ of dimension $1$ corresponds to the number of non-trivial
zeros for the quadratic equation $q_L=0$.
%$q_L(t_1,\ldots,t_{2\sigma})=0$. 
As above, set $q= t_1t_2+\ldots+t_{2\sigma-1}t_{2\sigma}$.
The number of zeros for $q=0$ 
is given by  $2^{2\sigma-1} +   2^{\sigma-1}$,
whereas $q+t_1^2+t_2^2=0$
%$q(t_1,\ldots,t_{2\sigma})+t_1^2+t_2^2=0$ 
has $2^{2\sigma-1} -  2^{\sigma-1}$ zeros (see for example \cite{Small 1991}, Theorem 4.16).
The assertion follows from \eqref{discriminat quadratic form}.  
\qed

\medskip
For   supersingular K3 surfaces $X$ in arbitrary characteristics $p>0$, the Picard lattice $L=\Pic(X)$ has rank $\rho=22$,
the intersection form is even, and the discriminant group has order $|L^*/L|=p^n$ for some $n\geq 0$.
Artin \cite{Artin 1974} showed   that the lattice $L$ is $p$-elementary, with even exponent $n=2\sigma$,
and least for odd characteristics. The required results on flat cohomology where established a little later by 
Milne \cite{Milne 1976}.
Rudakov and Shafarevich extended this to $p=2$ and established that the discriminant bilinear form $b_L$
is alternating (Rudakov and Shafarevich \cite{Rudakov; Safarevic 1979}, Theorem 3).
In turn, the   Picard lattice $L=\Pic(X)$ is determined up to isometry by 
the numbers  $r=22$ and   $1\leq\sigma\leq 10$. The latter is  called the \emph{Artin invariant}.
 
\begin{proposition}
\mylabel{artin invariant}
For the group scheme $G=\mu_2$,   the Artin invariant of the Kummer surface
$X=\Km(C\times C)$ is an integer  $1\leq \sigma\leq 3$, and each such number  
occurs.
\end{proposition}

\proof 
Consider the fibration $ X\ra \PP^1$ induced from  $\pr_2:C\times C\ra C$,
and write $K=k(\PP^1)$ for the function field of the projective line. We saw in Corollary \ref{rational points} that
the number of $K$-rational points in the generic fiber $X_K$ takes the form  $|X_K(K)|=2^{2+m}$ with $0\leq m\leq 2$.
Moreover,  each such integer can be realized with suitable vector fields $\delta\in H^0(C\times C,\Theta_{C\times C})$.

From Proposition \ref{distinguished curves} we deduce that  $ X\ra\PP^1$ has exactly five singular fibers, all of which 
are of type $\I_0^*$.
In turn, the trivial lattice $L=T(X/\PP^1)$ inside $P=\Pic(X)$    
has discriminant $\disc(L)=-2^{2\cdot 5}$, whereas the full Picard lattice has $\disc(P)=-2^{2\sigma}$. 
Finally, the index for the sublattice
is given by $[P:L]=2^{2+m}$. This gives $10=2\sigma+2(2+m)$, and thus $\sigma=3-m$.
The assertion now follows.
\qed

\medskip
Let $L\subset \Pic(X)$ be the sublattice generated by the distinguished curves in Figure \ref{curve configuration}.
This  lattice has invariants $r=22$ and $\sigma=3$, and occurs for all Kummer surfaces $X=\Km(C\times C)$
associated with group scheme $G=\mu_2$. By Lemma \ref{overlattices},
there are   $n_2=27$ overlattices $L\subset L'$ with   invariant $\sigma'=2$.
Moreover, for fixed $L'$ there are five further  overlattices $L'\subset L''$  with invariant $\sigma''=1$.
In total,  there are $n_1=\frac{5n_2}{3}=45$  such  overlattices $L\subset L''$.
We saw above that at least one     $L'$ and one   $L''$ appears as Picard groups  for  
Kummer surfaces $X=\Km(C\times C)$. 
We do not know which of them actually occur in this way.

%===========================================================
\section{Characterization with configurations of curves}
\mylabel{Configurations}

Throughout this section we work over an algebraically closed ground field $k$ of characteristic $p\geq 0$.
Over the complex numbers, the classical Kummer surfaces $X=\Km(A)$ with group $G=\{\pm 1\}$
attached to an abelian surface $A$ can be characterized
by the existence of sixteen disjoint $(-2)$-curves (Nikulin \cite{Nikulin 1975}).
Those coming from   jacobians of   genus-two curves are characterized by a
$(16_6)$-configuration of thirty-two $(-2)$-curves (Nikulin \cite{Nikulin 1974}).
Our main result 
establish a similar characterization for the Kummer surfaces $\Km(C\times C)$ associated with group scheme $G=\mu_2$.

\begin{theorem}
\mylabel{configuration and kummer}
Let $X$ be a K3 surface containing a configuration of thirty  $(-2)$-curves
with normal crossings and dual graph 
as in Figure \ref{curve configuration}. Then the characteristic must be $p=2$, and $X$ is isomorphic to a Kummer surface
$\Km(C\times C)$ associated with group scheme $G=\mu_2$.
\end{theorem}

\proof
We first construct a jacobian quasielliptic fibration $f:X\ra\PP^1$. Consider the following divisors:
\begin{equation}
\label{reducible fibers}
C_0 + 2E_0 + E_1 +E_2 + E_3\quadand 2C_i+ E_{i,1}+ \ldots + E_{i,4},\quad 1\leq i\leq 4.
\end{equation}
These are pairwise disjoint, and each forms a singular fiber of type $\I_0^*$.
Let $f:X\ra\PP^1$ be the resulting genus-one fibration. It is jacobian because the curve $C_1'$ provides a section.
The trivial lattice $T(X/\PP^1)$
has rank $r\geq 22$. It follows that the K3 surface $X$ has Picard number $\rho=22$ and 
the  Mordell--Weil group $\MW(X/\PP^1)$ is finite.
Consequently the above five divisors
are the  reducible fibers, and we have $p>0$.
The Picard lattice $\Pic(X)$ has discriminant $-p^{2\sigma}$.
The sublattice $L\subset\NS(X)$ generated by the irreducible curves appearing in \eqref{reducible fibers}, 
together with a section is isomorphic to $U\oplus D_4^{\oplus 5}$ which has even discriminant $-2^{10}$, and we conclude $p=2$.
Furthermore, the fibration $f:X\ra\PP^1$ is quasielliptic (Rudakov and Shafarevich \cite{Rudakov; Safarevic 1979}, Proposition on page 150).

By symmetry, the curves
\begin{equation}
\label{reducible fibers'}
C_0' + 2E_0 + E_1 + E_2 + E_3\quadand 2C_j'+ E_{1,j}+ \ldots + E_{4,j},\quad 1\leq j\leq 4 
\end{equation}
give another such fibration $f':X\ra\PP^1$.
Using the dual graph, we compute the intersection number $f^{-1}(\infty)\cdot f'^{-1}(\infty)$ between
the fibers as
$$
(2C_1+\sum_{i} E_{1,i})\cdot (2C_1'+\sum_{j} E_{j,1})  =
(2C_1+E_{11})\cdot (2C'_1+E_{11}) =2.
$$
In turn, the resulting morphism $(f,f'):X\ra\PP^1\times\PP^1$ is an alteration of degree two,
which means a proper surjection between integral scheme whose generic fiber has length two.
Let $X\ra Z\ra\PP^1\times\PP^1$ be its Stein factorization.
The morphism contracts precisely the irreducible curves $C\subset X$ that are vertical for both fibrations.
These curves correspond to the white vertices in the Figure \ref{curve configuration}, thus form an 
ADE-configuration of the type $16A_1+D_4$. In turn, $Z$ is a normal K3 surface
with Picard number $\rho=2$, coming with a finite flat morphism $Z\ra\PP^1\times\PP^1$.

Let $z_{ij}\in Z$ be the images of the exceptional curves $E_{ij}\subset X$,
and let $z\in Z$ be the image of $E_0\cup\ldots\cup E_4$. The complement   
of these seventeen singular points is the regular locus $U=\Reg(Z)$. We  now construct
an invertible sheaf on $U$ whose class in $\Pic(U)$ has order two, such that
the resulting $\mu_2$-torsor will lead to  the desired normal Kummer surface.
Let $g:Z\ra\PP^1$ be the morphism induced by the fibration $f:X\ra\PP^1$, and 
consider the images $F_i\subset Z$ of the  curves $C_i \subset X$.
Then $F_0,2F_1,\ldots,2F_4$ are schematic fibers for $g$, and it follows that the Weil divisor
$A=2F_0-(F_1+\ldots+F_4)$ has order two modulo principal divisors.
Fix an identification $\O_Z(-2A)=\O_Z$. 
For the reflexive rank-one sheaf $\shF=\O_Y(A)$,
we obtain a canonical map $\shF^\vee\otimes\shF^\vee\ra \O_Y$, which endows the coherent sheaf $\shA=\O_Z\oplus\shF^\vee$
with the structure of a  $\O_Z$-algebra graded by the group $\ZZ/2\ZZ$. Let $\epsilon:\tilde{Z}\ra Z$ be the resulting finite $Z$-scheme,
which is irreducible and Cohen--Macaulay.
The $\ZZ/2\ZZ$-grading on $\shA$ corresponds to an action of the group scheme $G=\mu_2$ on $\tilde{Z}$, 
with quotient $Z=\tilde{Z}/G$. Over the regular locus $U=\Reg(Z)$, the action is free,
and the quotient map becomes a $G$-torsor.

Let $g':Z\ra\PP^1$ be the morphism induced by the fibration $f':X\ra\PP^1$.
The situation is actually symmetric in $g$ and $g'$: Using the
dual graph given in Figure \ref{curve configuration}, one sees that the divisors 
$$
2C_0-(C_1+\ldots+C_4)\quadand 2C_0'-(C_1'+\ldots+C'_4)
$$
have
the same intersection numbers with all  curves occurring in the dual graph. Since these generate 
the Picard group up to finite index and $\Pic^\tau(X)=0$, we conclude that the above divisors differ by a prime divisor.
Let $F_i'\subset Z$ be the images of the curves $C_i'\subset X$, and set $A'=2F_0'-(F'_1+\ldots+F'_4)$.
Then $2F_0-(F_1+\ldots+F_4)$ and $2F_0'-(F_1'+\ldots+F'_4)$ also differ by a principal divisor.
In turn, our  double covering $\epsilon:\tilde{Z} \ra Z$ defined with $\shF=\O_Z(A)$
is equivariantly isomorphic to the double covering
defined with $\shF'=\O_Z(A')$.

\newcommand{\tZ}{{\tilde{Z}}}
The main task now is to identify $\tZ$ with the self-product $C\times C$ of the rational cuspidal curve.
We start by computing the Euler characteristic for the structure sheaf,
which boils down to compute $\chi(\shF)$ on the normal K3 surface $Z$.
First note that the half-fibers $F_i$, $1\leq i\leq 4$ are copies of the projective line:
To see this, write $h:X\ra Z$ for the contraction. For each $A_1$-singularity $z_{ij}\in Z$, the schematic fiber
is given by 
$h^{-1}(z_{ij})=C_{ij}$, according to \cite{Artin 1966}, Theorem 4, and this  ensures that  the induced morphism 
$h:C_i\ra F_i$ is an isomorphism.
Consider the disjoint
union $F=F_1\cup\ldots\cup F_4$.
The short exact sequence $0\ra\O_Z(-F)\ra\O_Z\ra\O_F\ra 0$ immediately
gives $\chi(\O_Z(-F))=\chi(\O_Z)-4\chi(\O_{\PP^1})= -2$.
In contrast to the half-fibers, the fiber $F_0$ is a a copy of the rational cuspidal curve:
The $D_4$-singularity $z\in Z$ has schematic fiber $h^{-1}(z)=2E_0+E_1+E_2+E_3$, again by \cite{Artin 1966}, Theorem 4.
In turn, $h^{-1}(z)\cap C_0= C_0\cap 2E_0$ is a local Artin scheme of length two, which is mapped to the closed
point $z\in F_0$ under the induced morphism $C_0\ra F_0$. It follows that $F_0$ is the rational cuspidal curve.
Since $F_0$ is a fiber, we furthermore have $2F_0=F_0\otimes k[\epsilon]$, where
$\epsilon$ is an indeterminate subject to $\epsilon^2=0$.
The short exact sequence $0\ra \O_Z(-F)\ra\O_Z(2F_0-F)\ra\O_{2F_0}\ra 0$
yields $\chi(\shF) = \chi(\O_Z(-F)) + \chi(\O_{2F_0}) = -2+2\cdot 0=-2$.
In turn, we get $\chi(\O_\tZ)=\chi(\O_Z) +\chi(\shF) = 2-2=0$.

Next  note that $\tZ$ is reduced: If not, the structure morphism $\epsilon:\tZ\ra Z$ admits
generically a section. By the Valuative Criterion for proper morphism, such a generic section
extends over an open subset $V\subset Z$ containing all points of codimension one.
Thus $\epsilon^{-1}(V)\ra V$ is a trivial $G$-torsor, and it follows that the invertible sheaf $\shF|U$
is trivial. In turn, the Weil divisor $A=2F_0-(F_1+\ldots+F_4)$ is principal over $V$, hence
on $Z$. However, this Weil divisor is not principal at each rational double point of type $A_1$,
contradiction. Summing up, our surface $\tZ$ is integral.

Furthermore,  we observe that the dualizing sheaf is isomorphic to the structure sheaf.
Indeed, we have $\omega_Z=\O_Z$, whence $\omega_{\tZ}$ is trivial over the open
set $\epsilon^{-1}(U)$. Since $\omega_Y$ is Cohen--Macaulay, we must have $\omega_{\tZ}=\O_\tZ$.
We now get the cohomological invariants: $h^0(\O_\tZ)=1$ because $\tZ$ is integral,
$h^2(\O_\tZ)=1$ by Serre duality, and finally $h^1(\O_\tZ)=2$ according to the computation of the
Euler characteristic.

\newcommand{\tg}{\tilde{g}}
Consider the composition $\tg:\tZ\ra\PP^1 $ of the double covering $\tZ\ra Z$ with
the fibration $g:Z\ra \PP^1$, and let $D=\Spec\tg_*(\O_{\tZ})$ be the Stein factorization.
Then $D$ is an integral curve, which turns out to be non-normal.
The morphism $D\ra\PP^1$ is a finite universal homeomorphism, and we claim that
it has degree two.
Indeed, the sheaf $\shF$ and the structure sheaf $\O_Z$ become isomorphic over the generic geometric fiber $S=f^{-1}(\bar{\eta})$,
and it follows that the $G$-torsor $\epsilon:\tZ\ra Z$ becomes trivial when pulled back to $S$.
It follows that the locally free sheaf $\tg_*(\O_{\tZ})$ has rank two.

We   observed at the beginning that  $f':X\ra \PP^1$ is quasielliptic.
So besides the five reducible fibers corresponding to  \eqref{reducible fibers'}, all closed fibers are copies of the rational cuspidal curve.
The corresponding fibers $Z_a=g'^{-1}(a)$, $a\in\PP^1$ on the normal K3 surface $Z$ are contained in $\Reg(Z)$,
and the induced torsor $\tZ_a=\epsilon^{-1}(Z_a)\ra Z_a$ is trivial. In turn,
the reduction $(\tZ_a)_\red$ is another copy of the rational cuspidal curve.
Taking degrees in the commutative diagram
$$
\begin{CD}
(\tZ_a)_\red	@>\epsilon>>	Z_a\\
@V\tg VV						@VVgV\\
D				@>>\can>			\PP^1,
\end{CD}
$$
we see that the finite dominant morphism $\tg:(\tZ_a)_\red\ra D$ is birational. In particular 
$h^1(\O_D)\geq 1$ holds.  

By symmetry, the above reasoning also applies  for the composition $\tg':\tZ\ra\PP^1$ of the double covering
 $\epsilon:\tZ\ra Z$ with
the other fibration $f':Z\ra\PP^1$, and the ensuing Stein factorization $D'=\Spec \tg'_*(\O_{\tZ})$.
Consider the resulting morphisms $\tZ\ra D$ and $\tZ\ra D'$ and the ensuing diagonal morphism
$\varphi:\tZ\ra D\times D'$ between integral schemes, which is proper and dominant.
We claim that it is birational: In the commutative diagram
$$
\begin{CD}
\tZ	@>\epsilon>>		Z\\
@V\varphi VV			@VV(g,g')V\\
D\times D'	@>>>	\PP^1\times\PP^1,
\end{CD}
$$
the upper   map has degree two, the right map has degree two,
and the lower map has degree four. If follows that $\deg(\varphi)=1$.

Next, we claim that the integral curves $D,D'$ have  $h^1(\O_D)=h^1(\O_{D'})=1$.
Seeking a contradiction, we assume that this  does not hold.
Without restriction, we have $h^1(\O_D)\geq 2$ and $h^1(\O_{D'})\geq 1$.
The canonical injection $H^1(\O_D)\subset H^1(\tZ,\O_\tZ)$  must be an equality, by dimension reasons.
To proceed, consider a fiber $\tZ_a=\tg^{-1}(a)$ such that the induced projection $g':\tZ_a\ra D'$
is birational. Then the composite map $H^1(D',\O_{D'})\ra H^1(\tZ,\O_\tZ)\ra H^1(\tZ_a,\O_{\tZ_a})$
is surjective. Hence there is a   cohomology class $\alpha\in H^1(\tZ,\O_\tZ)$
whose restriction to the fiber $\tZ_a=g^{-1}(a)$ is non-zero.
On the other hand, any cohomology class lies in the image of $g^*:H^1(D,\O_D)\ra H^1(\tZ,\O_\tZ)$,
whence vanishes on $Z_a$, contradiction.
Summing up, we have $h^1(\O_D)=h^1(\O_{D'})=1$.

Since the morphism $D\ra\PP^1$ and $D'\ra\PP^1$ are purely inseparable, we infer  that both $D,D'$ are copies of
the rational cuspidal curve. Summing up, we have a   birational morphism $\varphi:\tZ\ra C\times C$
between proper integral schemes, which are Gorenstein with trivial dualizing sheaves.
Consider the resulting conductor square
$$
\begin{CD}
R	@>>> 	\tZ\\
@VVV		@VVV\\
B	@>>>	C\times C,
\end{CD}
$$
where $B\subset C\times C$ is defined by the annihilator ideal for  
  $\varphi_*(\O_\tZ)/\O_{C\times C}$, and $R=\varphi^{-1}(C)$ is its schematic preimage.
If non-empty, the schemes $B,R$ are equidimensional of dimension one, and without embedded components,
because both $\tZ$ and $C\times C$ are Cohen--Macaulay.
The adjunction formula gives $\omega_\tZ = \varphi^*(\omega_{C\times C})(-R)$.
It follows that $R$ and hence also $B$ is empty, thus $\varphi:\tZ\ra C\times C$ is an isomorphism.

Summing up, our normal K3 surface $Z$ is isomorphic to the quotient
for a faithful action of $G=\mu_2$ on the self-product $\tZ=C\times C$.
In other words, the K3 surface $X$ is a Kummer K3 surface $\Km(C\times C)$ formed 
with the group scheme $G=\mu_2$.
\qed

%===========================================================
\section{Characterization with Artin invariants}
\mylabel{Artin invariant}

Let $k$ be an algebraically closed ground field of characteristic $p\geq 0$,
and let $X$ be a   K3 surface.
Then $\Pic(X)$ is a free abelian group of rank $\rho\leq 20$ or $\rho=22$ endowed with a non-degenerate
intersection pairing.

\begin{lemma}
\mylabel{picard isometry}
Suppose  $X$ and $X'$ are K3-surfaces such that there is an isometry between Picard groups.
Let $C_1,\ldots,C_r\subset X$ be  $(-2)$-curves.
Then there are $(-2)$-curves $C'_1,\ldots,C'_r\subset X'$ so that
the intersection matrices $N=(C_i\cdot C_j)$ and $N'=(C_i'\cdot C'_j)$ coincide.
\end{lemma}

\proof
By assumption there is an  isomorphism $f:\NS(X)\ra \NS(X')$
respecting the intersection pairing. According to Rudakov and Shafarevich \cite{Rudakov; Shafarevich 1983}, Section 3, Proposition  one may choose $f$
that the induced map between real vector space yields a bijection $\operatorname{Nef}(X)\ra \operatorname{Nef}(X)$
between the nef cones (see also Shimada and Zhang \cite{Shimada; Zhang 2007}, Proposition 3.1).
In turn, it     induces a bijection $\NE(X)\ra\NE(X')$ between  the cones of curves, which are dual to the nef cones.
But on any smooth projective surface $S$, the extremal rays $\RR_{\geq 0}c\subset \NE(S)$ 
with $c^2<0$ correspond to the integral curves $C\subset S$ with $C^2<0$,
according to Koll\'ar \cite{Kollar 1995}, Lemma 4.12.
On K3 surfaces, these are  exactly the $(-2)$-curves.
In turn, our chosen map $f:\NS(X)\ra \NS(X')$ induces a bijection between the classes of $(-2)$-curves
on $X$ and $X'$, and respects their intersection numbers.
\qed

\medskip
Now suppose that $X$ is a K3-surface with $\rho=22$. Then we are in characteristic $p>0$,
the K3-surface is supersingular, and the Picard group $\Pic(X)$ is determined up to isometry
by the Artin invariant $1\leq \sigma\leq 10$.
 
\begin{theorem}
\mylabel{kummer artin}
Let $X$ be a supersingular K3 surface in characteristic $2$ with Artin invariant $\sigma\leq 3$.
Then $X$ is isomorphic to a   Kummer  surface $\Km(C\times C)$ associated with 
group scheme $G=\mu_2$.
\end{theorem}

\proof
According to Proposition   \ref{artin invariant}, there is a Kummer surface $X'$ associated 
with group scheme $G=\mu_2$
whose Artin invariant  coincides with the given $\sigma$.
Note that $X'$ depends on the chosen embedding $h':G\ra\Aut_{C\times C}$.
According to Proposition \ref{distinguished curves} it 
contains thirty $(-2)$-curves with simple normal crossings and dual graph as in Figure \ref{curve configuration}.
By Lemma \ref{picard isometry} such a configuration also exists on $X$. According to  Theorem \ref{configuration and kummer} the 
K3 surface $X$ must be a Kummer surface associated with group scheme $G$,  for another embedding $h:G\ra\Aut_{C\times C}$.
\qed

%===========================================================
\section{Enriques surfaces and K3-like coverings}
\mylabel{K3-like coverings}

Keum \cite{Keum 1990} showed that every Kummer surface $X=\Km(A)$ over the field  $k=\CC$ of complex numbers admits a free action
of the group $H=\ZZ/2\ZZ$. 
Hence the  quotient $Y=X/H$ is an Enriques surface, and the Kummer surface
arises as its \emph{K3-covering}.
The result  was extend to characteristic $p\neq 2$ by Jang \cite{Jang 2015}.
If the abelian surface is a product $A=E\times E'$, the \emph{Lieberman involutions}
$$
E\times E'\lra E\times E',\quad (a,a')\longmapsto (a+\zeta,-a'+\zeta')
$$
induce  such  free $H$-actions on $X=\Km(E\times E')$, where $\zeta,\zeta'$ denote non-zero 2-division points.

We now establish an analogue for Kummer surfaces associated with group schemes. The case of Artin invariant $\sigma=1$
was already treated by the first author \cite{Kondo 2018}.

\begin{theorem}
\mylabel{k3-like covering}
Let $X$ be a supersingular K3 surface with Artin invariant $\sigma\leq 3$
over an algebraically closed ground field $k$ of characteristic $p=2$.
Then there is a contraction $X\ra X'$ of twelve $(-2)$-curves such that the normal K3 surface $X'$
is the K3-like covering of some simply-connected Enriques surface $Y$.
\end{theorem}

\proof
According to Theorem \ref{kummer artin}, our  K3 surface can be written as a  Kummer surface $X=\Km(C\times C)$ associated with group scheme $G=\mu_2$.
By Theorem \ref{configuration and kummer} it contains a  configuration of thirty $(-2)$-curves $E_{ij}$, $C_r$, $C_r'$,  $E_s$
with dual graph in Figure \ref{curve configuration}. The curves
\begin{gather*}
\label{I16 fiber}
E_{11} + E_{14} + E_{41} + E_{42} + E_{22} + E_{23} + E_{33} + E_{34} + \sum_{i=1}^4(C_i+C_i') 
\end{gather*}
form a singular fiber of type $\I_{16}$. Let $f:X\ra\PP^1$ be the resulting elliptic fibration.
The connected curve $E_0+\ldots+E_3$ is vertical with respect to this fibration,
hence belongs to some fiber, say $f^{-1}(\infty)$.
Let $s\geq 1$ be the number of irreducible components in this fiber.
Then the trivial lattice $T(X/\PP^1)$ has rank $r\geq 15+(s-1)+2=s+16$.
Since this lattice has rank at most $\rho=22$, we infer that $s\leq 6$.
It follows that $f^{-1}(\infty)$ has Kodaira symbol $\I_n^*$ with $n\leq 1$,   
 and that it  is the only
non-reduced fiber.

Let $X\ra X'$ be the contraction of the eight disjoint curves $C_1,\ldots,C_4,C_1',\ldots,C_4'$ together with $E_0+\ldots+E_3$.
Then $X'$ is a normal K3 surface with eight rational double points of type $A_1$, together with a rational double point
of type $D_4$. Moreover,   $f:X\ra\PP^1$ descends to a fibration $f':X'\ra\PP^1$.
We will apply a result of the second author (\cite{Schroeer 2017}, Theorem 6.4) to see that $X'$ is a K3-like covering;
for this we have to verify certain conditions (E0), (E2), (E5) and (E6):

First note that the fibers of $f':X'\ra \PP^1$ are reduced, because $f^{-1}(\infty)$ is the only non-reduced
fiber for $f:X\ra \PP^1$ and all components  with multiplicity $m\neq 1$ are contracted by $X\ra X'$.
Thus condition (E0) is satisfied.

Recall that the \emph{Tjurina number} of an isolated  hypersurface singularity given by a  power series equation
$f(t_1,\ldots,t_n)=0$
is the colength of the ideal $\ideala\subset k[[t_1,\ldots,t_n]]$ generated by $f$ and its partial derivatives
$\partial f/\partial t_i$, compare the discussion in \cite{Schroeer 2017}, Section 1.
Each   $A_1$-singularity is formally isomorphic
to $z^2+xy=0$, whence its Tjurina number is $\tau=2$.
According to Artin's computations, there are exactly two formal isomorphism classes of $D_4$-singularities,   one
of which is simply-connected
(\cite{Artin 1977}, Section 3). 
In light of the universal homeomorphism $\PP^1\times\PP^1\ra (C\times C)/G$,
our rational double point must be simply-connected, hence is formally given by 
$z^2+x^2y+xy^2=0$,  with Tjurina number $\tau =8$.  
We see that the Tjurina numbers of the singularities add up to $n=24$, so (E2) holds. Moreover,  all singular local rings $\O_{X',a}$  
are Zariski singularities, hence (E6) is true. In particular the stalks $\Theta_{X',a}$
of the tangent sheaf are isomorphic to $\O_{X',a}^{\oplus 2}$,
see \cite{Schroeer 2017}, Corollary 1.6.

The configuration $C_1+E_{13}+C'_3+E_{43}+C_4+E_{42}+C'_2+E_{12}$ defines  another elliptic fibration
$g:X\ra\PP^1$, which descends to an elliptic fibration $g':X'\ra \PP^1$.
It follows that there is an elliptic curve $F\subset X'$ that is horizontal for   $f':X'\ra \PP^1$.
This establishes condition (E5).

Summing up, \cite{Schroeer 2017}, Theorem 6.4 applies, and we conclude that
$\Theta_{X'/k}=\O_{X'}^{\oplus 2}$, all members
of the restricted Lie algebra $\lieg=H^0(X,\Theta_{X'/k})$ are $p$-closed, and for almost all 
non-zero vectors $\delta \in \lieg$ the corresponding local group scheme $H\subset\Aut_{X'/k}$ of order $p$
acts freely, such that the quotient $Y=X'/H$ is an Enriques surface, with $X'$ as the K3-like covering.
\qed

%===========================================================

\end{document}